\documentclass[12pt,centertags,oneside]{amsart}
\usepackage{amsmath,amstext,amsthm,amscd,typearea,hyperref}
\usepackage{amssymb}
\usepackage{a4wide}
\usepackage[mathscr]{eucal}
\usepackage{mathrsfs}
\usepackage{typearea}
\usepackage{charter}
\usepackage{pdfsync}
\usepackage[a4paper,width=16.2cm,top=3cm,bottom=3cm]{geometry}

\numberwithin{equation}{section}

\newtheorem{theorem}{Theorem}[section]

\newtheorem{lemma}[theorem]{Lemma}

\def\IH{{\mathbb H}} 
 
\def\IR{{\mathbb R}} 
\def\IQ{{\mathbb Q}}

\def\IZ{{\mathbb Z}} 
\def\IN{{\mathbb N}} 
 
\def\tr{{\rm trace}}

\def\G{\Gamma}

\usepackage{graphicx}

\emergencystretch=\maxdimen
\title[]{On Thin Heckoid and Generalised Triangle Groups in $PSL(2,\mathbb{C})$}

\author{Alex Elzenaar, Gaven Martin and Jeroen Schillewaert}
\address{Institute for Advanced Study, Massey University, New Zealand}
\email{g.j.martin@massey.ac.nz}

\date{}
\begin{document}

\maketitle

\begin{abstract} We provide a brief overview of our upcoming work identifying all the thin Heckoid groups in $PSL(2,\mathbb{C})$. Here we give a complete list of the $55$ thin generalised triangle groups of slope $1/2$. This work was presented at the conference Computational Aspects of Thin Groups, IMSS, Singapore  and presents an application of joint work initiated with Colin Maclachlan
\end{abstract}

\medskip

\noindent {\bf Classification AMS 2020}: 20H10 30F40 57S30

\medskip



\medskip

Let $\Gamma < GL_n(\IZ)$ be a finitely generated subgroup and $G = Zcl(\Gamma)$ its Zariski closure. $\Gamma$ is a  {\it thin group} if the index of $\Gamma$ in the integer points $G(\IZ)$ is infinite.   
For a subgroup of $PSL(2,\mathbb{C})\cong {\rm Isom}^+(\IH^3)$ this definition simplifies considerably.

\begin{lemma}[{Definition}] Let $\Gamma\subset PSL(2,\mathbb{C})$ be a thin group.  Then $\Gamma$ is an infinite index subgroup of an arithmetic Kleinian group.  Thus the invariant trace field 
\[ k\Gamma^{(2)} = \IQ(\{\tr^2(h):h\in \Gamma\})\]
has one complex place and the invariant quaternion algebra is ramified at real places.
\end{lemma}

Let $\Gamma\subset PSL(2,\mathbb{C})$, discrete and not virtually abelian. $\Gamma$ is a {\it Heckoid group} if $\Gamma$ can be generated by two elements of finite order (or parabolic) and is not freely generated by these elements. 
\[ \Gamma \cong \langle a,b:a^p=b^q=w_1(a,b)=\cdots=w_n(a,b) \rangle \]
Classical examples are the triangle groups
$\Delta(p,q,r)=\langle a,b:a^p=b^q=(ab)^r=1 \rangle$, $\frac{1}{p}+\frac{1}{q}+\frac{1}{r}<1$, 
 which embed as lattices in $PSL(2,\IR)$.  A {\it generalised triangle group} is a Heckoid group with presentation of the form
$\langle a,b:a^p=b^q=w(a,b)^n=1 \rangle$, $n\geq 2.$
 
 The main result we present is a complete identification of the thin Heckoid subgroups of ${\rm Isom}^+(\IH^3)$.
 
 \begin{theorem} There are finitely many thin Heckoid subgroups of ${\rm Isom}^+(\IH^3)$.  The degree of the invariant trace field is never more than $8$. The order of the torsion is never more than 30. Both bounds are sharp.
 \end{theorem}

In all there are about 150 such groups (up to Nielsen equivalence and conjugacy).  Precise tables listing each example and the associated arithmetic data will appear.

\section{A motivational example.}\label{parasec}

In this section we carefully go through the simplest example,  that of Heckoid groups generated by two parabolic elements. It is simplest because the existence of parabolic elements implies that the degree of the invariant trace field is at most two. Bounding the degree of this trace field in the general case is really one of the most difficult issues.    This is a good example in the sense that the general strategy for finding all the groups in question follows this argument and it identifies the key steps quite clearly.
  
 \subsection{The parabolic Heckoid groups.} 
In important  papers the Japanese groups Aimi, Lee, Sakai, and M. Sakuma, and also H. Akiyoshi, K. Ohshika, J. Parker, M. Sakuma and H. Yoshida established a conjecture of Agol's on groups generated by {\it two parabolic elements}; they are near relatives of two bridge knots and links, see \cite{ALSS} and \cite{AOPSY}.

The set of all discrete Kleinian groups freely generated by two parabolic elements is the closure of the Riley slice -- the unbounded region laminated by Keen--Series pleating rays in Figure 1.  This region has been studied for decades and a lot is known about it,  see  \cite{EMS1,KS,KS2} and the references therein.  Our paper \cite{EMS2} has more historical commentary. The aforementioned work \cite{ALSS,AOPSY} completely describes the discrete groups in the complement of the Riley slice. 

We use all of this work to give a proof of the following theorem outlining the main steps.  This arithmetic result appears in \cite{GMM}.

\begin{theorem} Let $\Gamma$ be a Kleinian group generated by two parabolics which does not split as a free product $\Gamma \not\cong \IZ*\IZ$. Among this infinite family of groups
\begin{itemize} 
\item Exactly $4$ are arithmetic : one knot and three links. 4,5 \& 6 crossings. Torsion free!\\
\vspace{-0.2cm}
\begin{center}
\scalebox{0.3}{\includegraphics{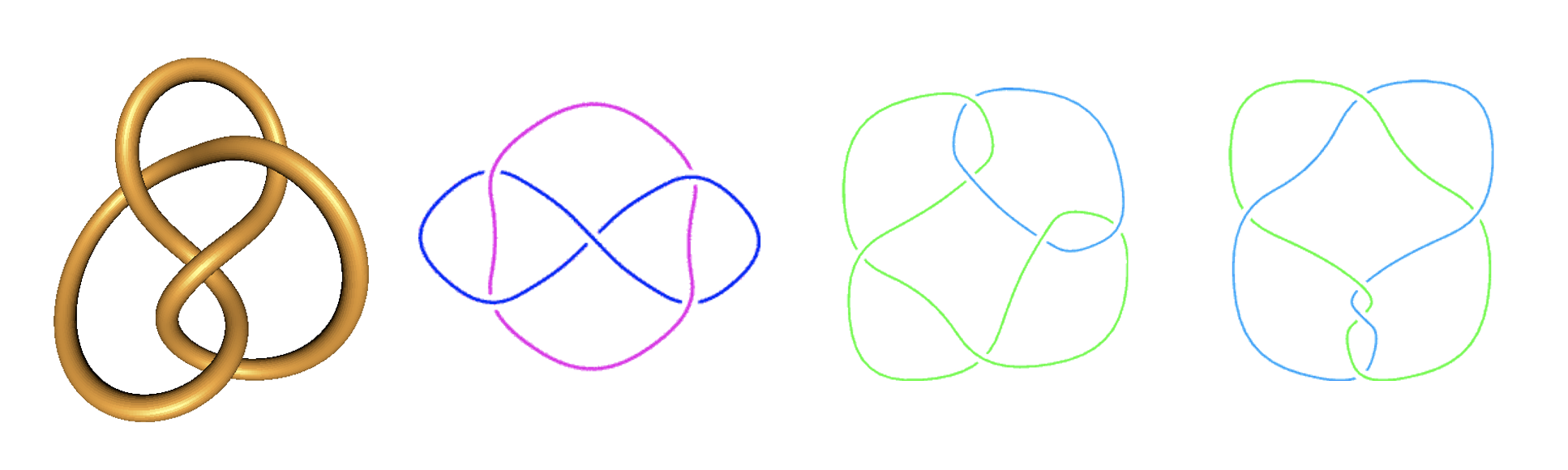}}
\end{center}
\item Exactly $3$ are thin,  have torsion and one or two hyperfinite vertices.
\vspace{-0.2cm}
\begin{center}
\scalebox{0.2}{\includegraphics{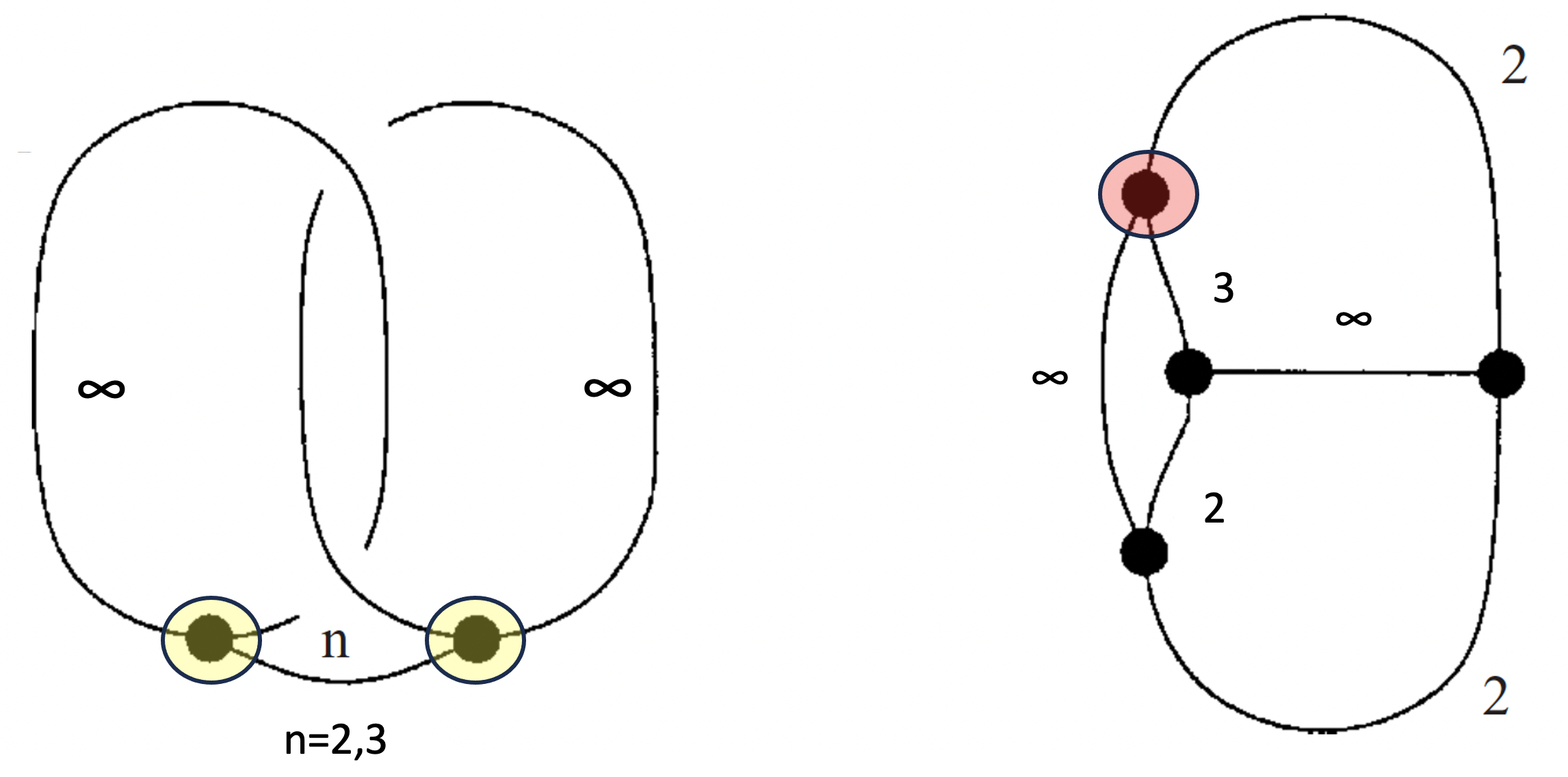}}
\end{center}
\end{itemize} 
\end{theorem}

\noindent {\bf Proof.}
\begin{enumerate}
\item Representations up to conjugacy:
\[ \Gamma \cong \Big\langle \left(
\begin{array}{cc}
1  & 1 \\
 0 & 1 \\
\end{array}
\right), \left(
\begin{array}{cc}
1  & 0 \\
 \rho  & 1 \\
\end{array}
\right)\Big\rangle\]
The single conjugacy invariant is $\gamma = \tr [A,B]-2 = \rho^2$.
\item Subgroups of a nonuniform arithmetic group, so the invariant trace field $k\Gamma=\IQ(\tr^2(h):h\in \Gamma)$
is quadratic, and $k\Gamma = \IQ(\gamma)$.
\item The group does not split and so lies outside the quasiconformal deformation space of $\IZ_\infty*\IZ_\infty$. That is $ \gamma \not\in {\mathcal  R}$,  the {\it Riley slice}.
\end{enumerate}
Thus we seek all complex quadratic integers in the region bounded by the Riley slice.  These are easy to completely enumerate.
\begin{center}
\scalebox{0.35}{\includegraphics{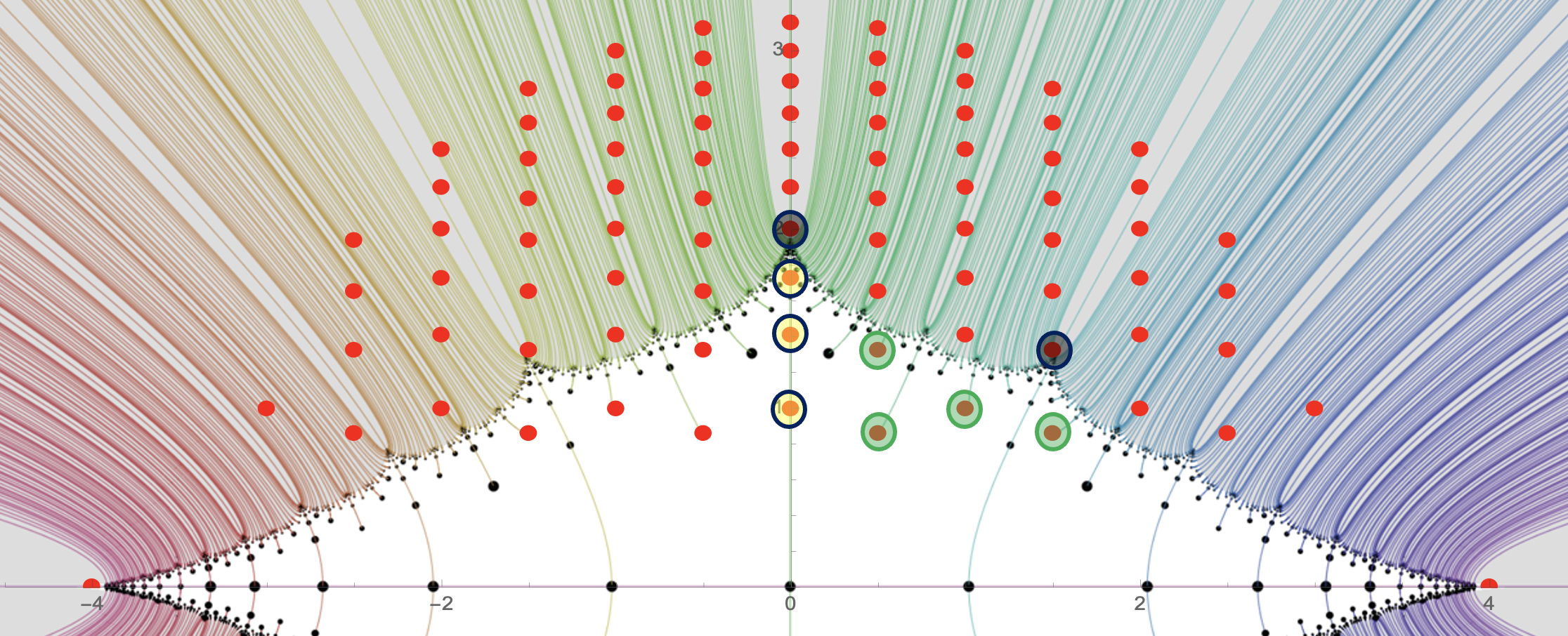}}
\end{center}
\medskip
\noindent{\bf Figure 1.} {\it The Riley slice is the foliated unbounded region. The enumerated points are illustrated including all the points in the complement of the Riley slice.  Image of Riley slice courtesy of Yasushi Yamashita.}

\medskip
There is a natural symmetry of the Riley slice $\rho\leftrightarrow -\rho$ which obviously yields the same groups. Reducing by this symmetry  
we see(!) that there are $7$ points in this set we must consider and answer the following questions:
\begin{itemize}
\item Are these representations of discrete groups ? 
\item Subgroups of arithmetic groups ?
\item Can we identify them algebraically ?
\end{itemize}
In the case at hand this is fairly straightforward since there are only $7$ points and the trace field is a discrete ring.  The results we found are presented.  However there are some issues.  Two of the points ($\frac{1}{2}(1+i\sqrt{7})$ and $i\sqrt{2}$) are precisely on the boundary of the Riley slice - a quite  complicated fractal. These {\it cusp groups} at the end of pleating rays are readily identifiable and are discrete with limit set $\Lambda(\Gamma)$ a circle packing and quotient $\big(\hat{\mathbb{C}}\setminus \Lambda(\Gamma)\big)/\Gamma$ a pair of triply punctured spheres.  In general it is impossible to decide if a point is on the boundary,  unless we are lucky enough to identify it as a cusp group.  Fortunately this always happens and is a consequence of our work with Chesebro and Schillewaert on the algebraic structure of Heckoid groups building on \cite{ALSS,AOPSY,orbifold}.  Here we only have to {\it believe} this to be the case,  and search for a Farey word relator.

\section{The general case} Following this proof raises the following issues when when we try this for all $p$ and $q$. 

\subsection{Quasiconformal deformation spaces.} We must be able to ``computationally'' identify the quasiconformal deformation space of $\IZ_p*\IZ_q$, denoted ${\mathcal  R}^{p,q}$,  and \underline{effectively} bound it to limit searches. This is the subject of \cite{EMS4,EGMS} generalising the Keen--Series theory to the case of groups with two generators of finite order,  and earlier results bounding the space of discrete and free groups -- in particular Lyndon--Ullman's work.

\subsection{Arithmetic restrictions on the commutator.}  After choosing suitable matrix generators $A,B$ of primitive elliptics (minimal rotation angle) for a given pair of orders $p$ and $q$ there is a single conjugacy invariant identifying the group \cite{GM1},  this is
\[ \gamma = \tr[f,g]-2\]  
Arithmetic restrictions on $\gamma$ to imply $\Gamma=\langle A,B\rangle$ is a subgroup of an arithmetic Kleinian group were first developed in \cite{GMMR} and refined to this setting in \cite{MM,MMpq}.  These are encapsulated in the following theorem.
\begin{theorem}\label{MMthm}
Let $\G = \langle  f , g \rangle $ be a Kleinian group with $f$ of order $p$ and $g$ of order $q$, $3\leq q\leq p$.  Let $\gamma(f,g) = \gamma \in \mathbb{C} \setminus \IR$.  Then $\Gamma$ is an arithmetic Kleinian group if and only if 
\begin{enumerate}
\item $\gamma$ is an algebraic integer,
\item $\IQ(\gamma) \supset L = \IQ(\cos 2 \pi/p, \cos 2 \pi/q)$ and $\IQ(\gamma)$ is a number field with exactly one complex place,
\item if $\tau : \IQ(\gamma) \rightarrow \IR$ such that $\tau |_L = \sigma$, then
\[
 { - \sigma(4\sin^2\frac{\pi}{p}\sin^2\frac{\pi}{q}) < \tau(\gamma) < 0},
\]
\item From Fricke identity: the quadratic polynomial 
{  \[ x^2 - 16\cos^2\frac{\pi}{p}\cos^2\frac{\pi}{q} \, x 
  +  16\cos^2\frac{\pi}{p}\cos^2\frac{\pi}{q}(4\sin^2\frac{\pi}{p} -4\cos^2\frac{\pi}{q}  -\gamma) = 0 \]}
  factorises over $\IQ(\gamma)$,
\item $\Gamma$ has finite co-volume (without this,  we get that $\Gamma$ is a subgroup of an arithmetic group).
\end{enumerate}
\end{theorem}

\subsection{Finding possible candidates.}  In \cite{MM} it is shown that in any compact region of the plane,  there are only finitely many triples $(p,q,\gamma)\in \IN_{\geq 2} \times \IN_{\geq 3} \times \mathbb{C}$ which satisfy the conclusions of Theorem \ref{MMthm}.  In the special case of thin groups an analysis of the hyperfinite vertex stabilisers shows them to be arithmetic triangle groups,  completely classified by Takeuchi \cite{Tak} and hence we get restrictions on $p,q$ not available in the non-thin case. However,  even with sharp bounds on the size of the deformation space anything with degree more than $10$ is computationally infeasible on a reasonable desktop machine running Mathematica, as we are using. The central issue is using Theorem \ref{MMthm} to obtain good degree bounds on $\IQ(\gamma)$. If $\gamma$ is complex, then $k\Gamma^{(2)}=\IQ(\gamma)$,  while if $\Gamma$ is real it falls into the class of groups  discussed in \cite{MM2,MM3} and we reanalyse those cases seeking thinness (as opposed to arithmeticity which is easier).

In our example of \S \ref{parasec} we had a priori that this degree was two.  In general we use a number of computational techniques building on our earlier work \cite{MMpq}.  In response to a question of ours  Flammang and Rhin \cite{FR} essentially worked through the case of $p=2$, $q=3$.   Different approaches are needed depending on the ``shape'' of the deformation space and how large $\gamma$ is.  We typically bound the discriminant and turn this into a degree bound using computation number theoretic databases such as \cite{JR}.  See e.g.  \cite{MSY} when one generator has order $4$.  

Once degree bounds are at hand we search for  possible monic integral polynomials in the most obvious way; obtaining bounds on the coefficients as elementary symmetric functions of the roots which we have good bounds for from Theorem \ref{MMthm}.  With more knowledge (such as a discriminant bound),  number theoretic resources such as databases of fields and discriminants \cite{JR} can be very useful.

\subsection{Discrete but not free.} A candidate value for $\gamma$ satisfying the conclusions of Theorem \ref{MMthm} definitely gives us a discrete group. Before the rather precise descriptions of ${\mathcal  R}^{p,q}$ mentioned above,  other ad hoc techniques were used to prove groups were free or not of finite covolume,  \cite{CMMO, Cooper, Zhang}.  In \cite{EMS4} we prove certain neighbourhoods of pleating rays lie completely within the deformation space and are tangent to its boundary meeting  at a cusp group. This gives an algorithm (really a process) to decide if a value $\gamma\in {\mathcal  R}^{p,q}$ as these neighbourhoods exhaust  ${\mathcal  R}^{p,q}$.  In our simple example we already found candidates on  the fractal boundary and we have not eliminated the possibility that a candidate (arising from purely arithmetic data) is a geometrically infinite boundary point.  Though our recent work establishes that this does not happen.

\subsection{Identify the group.} At this point we have at hand a group which is discrete but not free and we wish to identify it and decide if it is thin.   We {\it believe} that the group has a power of a  Farey word \cite{EMS3} as a relator,  and in the case of thin groups this must be a proper power $n\geq 2$.  Farey words are determined by a ``rational slope'' $\frac{r}{s}$ related to the slope of a two-bridge knot or link \cite{BZ} (the Farey word is roughly the Wirtinger relator) and thus we have a symbol (putatively) determining the group $(p,q;r/s,n)_i$.  The subscript $i$ separates  the groups into two classes,  the generalised triangle groups and the {\it pure} Heckoid groups.  Our analysis of the traces of Farey words, called {\it Farey polynomials} $F_{r/t}^{p,q}(z)$ \cite{EMS3} allows for a quick search through slopes (Farey sequences) with denominators not exceeding say $N$,  so possible relators of word length 2N, and there are about $\frac{3}{\pi} N^2$ of these (as opposed to the potential exponential growth of possibilities). Given our candidate value $\gamma$ we search for $r/s$ and $n$ so that
\[F_{r/t}^{p,q}(\gamma) =  2\cos \frac{k\pi}{n}, \quad (k,n)=1, \quad \gamma=\gamma(f,g).\]
This is rather easier than it might appear since discreteness implies this is equivalent to $F_{r/t}^{p,q}(\gamma)\in[-2,2]$.  It turns out (and is now proved) that $k\in \{\pm1,\pm2\}$ ($n$ odd in the second case).  Finally an ``Euler" formula from the symbol determines whether the group is thin or arithmetic. A presentation can also be determined (though this is quite a bit more complicated for the groups which are not generalised triangle groups.)

\section{$p=q=3.$} We illustrate the implementation of this process in the following sequence of figures.
\begin{center}
{\scalebox{0.2}{\includegraphics{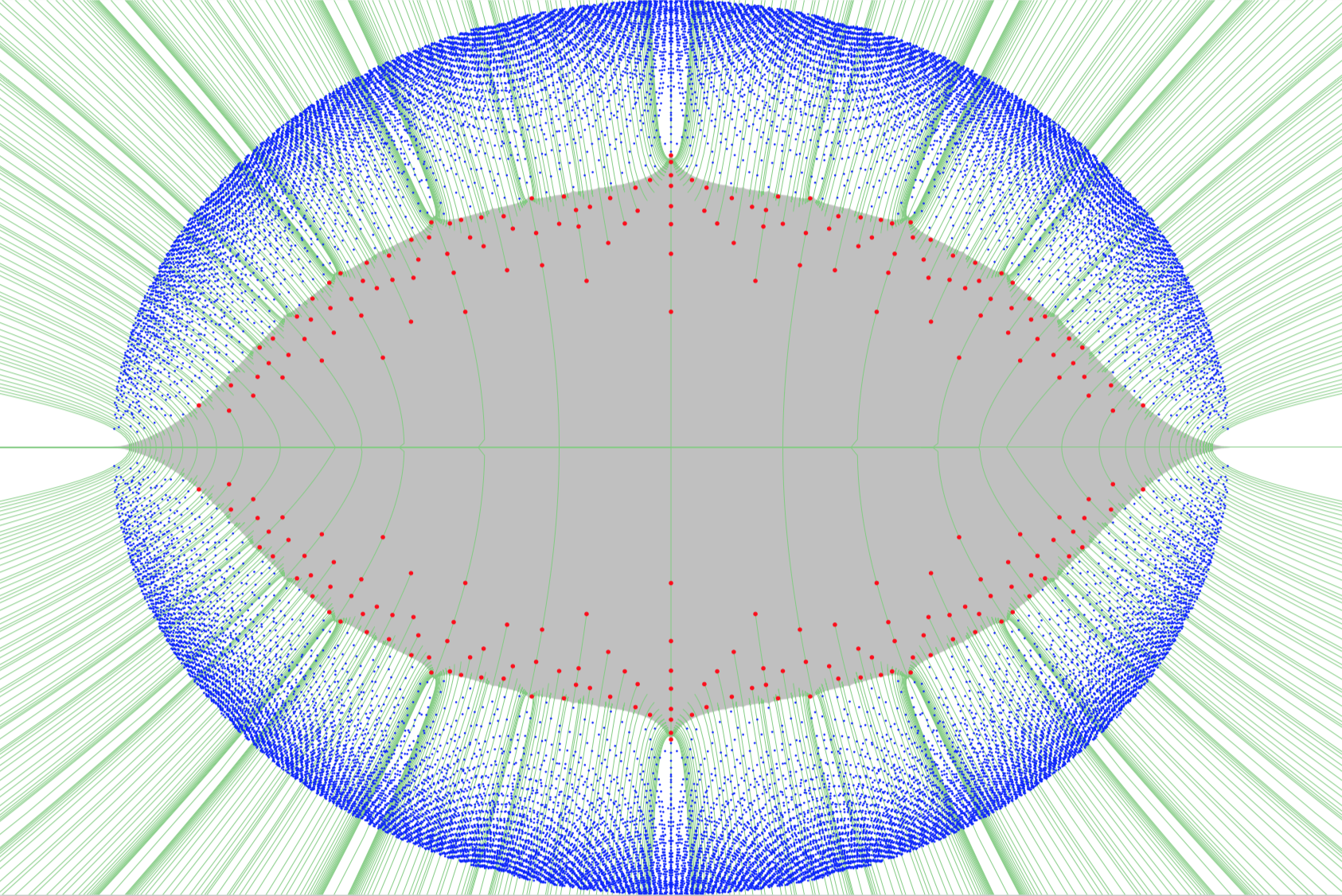}}}
\end{center}
\noindent{\bf Figure 2.} {\it Commutator values of groups generated by two elements of order three. 
\begin{itemize}
\item Shaded grey area is the {\em complement} of the quasiconformal deformation space of $\IZ_3*\IZ_3$. This deformation space is foliated by pleating rays.
\item All $15,909$ points satisfying the conclusions of Theorem \ref{MMthm} and are not obviously discrete and free via ping-pong type argument based on isometric circles (defining an ellipse). Total degree of $\IQ(\gamma)$ is no more than $10$. Some very close to the boundary.
\item $127$ red points are  either arithmetic or thin and Heckoid.
\item Thinness implies additional torsion - a triangle group representing a hyperfinite vertex. We use this to refine the candidate list.
\begin{center}
\scalebox{0.3}{\includegraphics{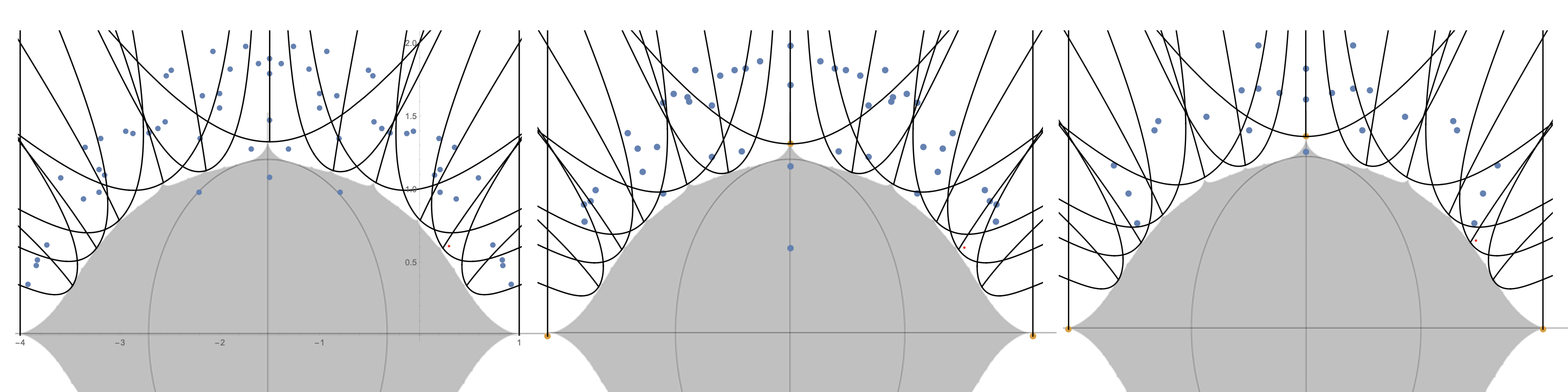}}
\end{center}
\noindent{\bf Figure 3.} From left $\sqrt{2}\in k\Gamma^{(2)}$, $\sqrt{3}\in k\Gamma^{(2)}$, $\sqrt{5}\in k\Gamma^{(2)}$. This refines the point sets which are then captured by pleating ray neighbourhoods (not all shown).
\end{itemize}

\section{An example: the thin generalised triangle groups of slope $\frac{1}{2}$.} 

As a specific example of these techniques we pressent a complete list of all thin generalised triangle groups of the type considered by several authors,  see  \cite{HMR,HMV,W} and also \cite{MM2}.  Those of slope $\frac{1}{2}$ all have orbifold type as illustrated in the statement of  Theorem 1.2,  lower left.

\medskip

We follow the following process in this special case.

\begin{itemize}  
\item We check that one of the vertices $(p,p,n)$ or $(q,q,n)$ is hyperfinite,  so the group cannot be a finite index subgroup of an arithmetic group.  
\item That the slope is $\frac{1}{2}$ implies the commutator  $[f,g]$ is elliptic,  and from \cite{MM3} we obtain
\item $\gamma=-2-2\cos\frac{\pi}{n}$,  thus from the arithmetic criterion Theorem 2.1 (4) we see  
\[ \alpha = 8 \cos \frac{\pi }{p} \cos \frac{\pi }{q}  \sqrt{4 \sin ^2 \frac{\pi }{p}  \sin ^2 \frac{\pi }{q}-2-2\cos\frac{\pi}{n} } \]  
\item We remove high degree candidates by considering Galois images of the ``discriminant'' term
\[ 4 \sin ^2 \frac{\pi }{p}  \sin ^2 \frac{\pi }{q}-2-2\cos\frac{\pi}{n} \]
\item E.g. if $n\in\{5,7,9,15\}$ and $p,q\neq n$ we must have
\[ 4 \sin ^2 \frac{\pi }{p}  \sin ^2 \frac{\pi }{q}-2-2\cos\frac{2\pi}{n} > 0\]
and obvious modifications if one or both of $p,q=n$ 
\item For each triple identified find the minimal polynomial for $\alpha$. 
\item Check the roots of this polynomial to see that there is exactly one complex conjugate pair.
\end{itemize}
 
 The $55$ thin generalised triangle groups of slope $\frac{1}{2}$ groups are the following: 
\begin{center} 
{\footnotesize
\begin{tabular}{|c|c|c|c|} 
\hline
1 & $(3,5;1/2,2)_1 $ &-275& $-11+9 z^2+z^4$\\ \hline
2 & $(3,6;1/2,2)_1 $ & -15& $15+z^2$\\ \hline
3 & $(3,8;1/2,2)_1 $ & -1792& $-28+20 z^2+z^4$\\ \hline
4 & $(3,12;1/2,2)_1 $ & -3312& $-23+26 z^2+z^4$\\ \hline
5 & $(4,6;1/2,2)_1 $ & -4& $36+z^2$\\ \hline
6 & $(5,5;1/2,2)_1 $ & -775& $-31+41 z^2+z^4$\\ \hline
7 & $(5,10;1/2,2)_1 $ & -1375& $-55+70 z^2+z^4$\\ \hline
8 & $(6,6;1/2,2)_1 $ & -7& $63+z^2$\\ \hline
9 & $(8,8;1/2,2)_1 $ & -448& $-112+88 z^2+z^4$\\ \hline
10 & $(12,12;1/2,2)_1 $ & -6768& $-47+110 z^2+z^4$\\ \hline  \hline
11 &$(3,4;1/2,3)_1 $ & -3& $12+z^2$\\ \hline
12 &$(3,6;1/2,3)_1 $ & -3& $27+z^2$\\ \hline
13 &$(4,4;1/2,3)_1 $ & -8& $32+z^2$\\ \hline
14 &$(4,6;1/2,3)_1 $ & -15& $60+z^2$\\ \hline
15 &$(5,5;1/2,3)_1 $ & -275& $-11+69 z^2+z^4$\\ \hline
16 &$(6,6;1/2,3)_1 $ & -11& $99+z^2$\\ \hline
17 & $(12,12;1/2,3)_1 $ & -3312& $-23+166 z^2+z^4$\\ \hline \hline
18 & $(3,3;1/2,4)_1 $ & -1984& $-31-2 z^2+z^4$\\ \hline
19 & $(3,4;1/2,4)_1 $ & -448& $-112+8 z^2+z^4$\\ \hline
20 & $(3,6;1/2,4)_1 $ & -4032& $-63+30 z^2+z^4$\\ \hline
21 & $(3,8;1/2,4)_1 $ & -12032& $-188+36 z^2+z^4$\\ \hline
22 & $(4,4;1/2,4)_1 $ & -1024& $-256+32 z^2+z^4$\\ \hline
23 & $(4,8;1/2,4)_1 $ & -1792& $-448+80 z^2+z^4$\\ \hline
24 & $(6,8;1/2,4)_1 $ & -16128& $-252+132 z^2+z^4$\\ \hline
25 & $(8,8;1/2,4)_1 $ & -7936& $-496+152 z^2+z^4$\\ \hline
\end{tabular}}
\end{center}
 
\begin{center} 
{\footnotesize
\begin{tabular}{|c|c|c|c|} 
\hline
26 &$(3,3;1/2,5)_1 $ & -475& $-19+2 z^2+z^4$\\ \hline
27 &$(3,4;1/2,5)_1 $ & -400& $-16+16 z^2+z^4$\\ \hline
28 &$(3,5;1/2,5)_1 $ & -1375& $-55+25 z^2+z^4$\\ \hline
29 &$(3,10;1/2,5)_1 $ & -3875& $-155+45 z^2+z^4$\\ \hline
30 &$(4,5;1/2,5)_1 $ & -2000& $-80+60 z^2+z^4$\\ \hline
31 &$(5,5;1/2,5)_1 $ & -2375& $-95+85 z^2+z^4$\\ \hline
32 &$(5,10;1/2,5)_1 $ & -275& $-275+130 z^2+z^4$\\ \hline
33 & $(10,10;1/2,5)_1 $ & -475& $-475+185 z^2+z^4$\\ \hline \hline
34 &$(3,3;1/2,6)_1 $ & -6768& $-47-2 z^2+z^4$\\ \hline
35 &$(3,4;1/2,6)_1 $ & -6336& $-176+8 z^2+z^4$\\ \hline
36 &$(3,6;1/2,6)_1 $ & -3312& $-207+30 z^2+z^4$\\ \hline
37 &$(3,12;1/2,6)_1 $ & -20592& $-143+50 z^2+z^4$\\ \hline
38 &$(4,4;1/2,6)_1 $ & -4608& $-512+32 z^2+z^4$\\ \hline
39 &$(4,6;1/2,6)_1 $ & -1728& $-432+72 z^2+z^4$\\ \hline
40 &$(4,12;1/2,6)_1 $ & -3312& $-368+104 z^2+z^4$\\ \hline
41 &$(6,12;1/2,6)_1 $ & -5616& $-351+162 z^2+z^4$\\ \hline
42 &$(12,12;1/2,6)_1 $ & -27504& $-191+206 z^2+z^4$\\ \hline
43 & $(24,24;1/2,6)_1 $ & -249495552& $-47+19060 z^2-4234 z^4+212 z^6+z^8$\\ \hline 
\end{tabular}}
\end{center}
 
\begin{center} 
{\footnotesize
\begin{tabular}{|c|c|c|c|} 
\hline
44 &$(3,7;1/2,7)_1 $ & -218491& $91-294 z^2+35 z^4+z^6$\\ \hline
45 &$(7,7;1/2,7)_1 $ & -487403& $203-1134 z^2+147 z^4+z^6$\\ \hline \hline
46 &$(3,3;1/2,9)_1 $ & -728271& $111-45 z^2-3 z^4+z^6$\\ \hline
47 &$(3,9;1/2,9)_1 $ & -1436859& $219-387 z^2+42 z^4+z^6$\\ \hline
48 &$(4,9;1/2,9)_1 $ & -1259712& $192-720 z^2+96 z^4+z^6$\\ \hline
49 &$(9,9;1/2,9)_1 $ & -728271& $111-927 z^2+186 z^4+z^6$\\ \hline 
50 &$(9,18;1/2,9)_1 $ & -334611& $459-1026 z^2+207 z^4+z^6$\\ \hline \hline
51 &$(5,5;1/2,10)_1 $ & -316000000& $-79+1458 z^2-1101 z^4+82 z^6+z^8$\\ \hline \hline
52 &$(3,3;1/2,15)_1 $ & -74671875& $-59+368 z^2-46 z^4-8 z^6+z^8$\\ \hline
53 &$(3,5;1/2,15)_1 $ & -454359375& $-359+638 z^2-241 z^4+22 z^6+z^8$\\ \hline
54 &$(3,15;1/2,15)_1 $ & -302484375& $-239+1888 z^2-646 z^4+47 z^6+z^8$\\ \hline
55 &$(5,15;1/2,15)_1 $ & -74671875& $-59+703 z^2-646 z^4+152 z^6+z^8$\\ \hline \hline
 \end{tabular}}
\end{center}


\end{document}